\documentclass{artjlt}

\usepackage{amsmath,amsfonts,amssymb,DLdef1-small}
\newcommand{\twol}[2]{\begin{tabular}{@{}c@{}}#1\\#2\end{tabular}}
\let\url\relax

\usepackage{hyperref}

\def\Subsection#1{\par\medskip\noindent%
\bf #1\quad\rm\ignorespaces}

\title{On integrations and cross ratios on
supermanifolds}                                     
\author{Dimitry Leites}                 
\lastname{Leites}  



\address{%
Dimitry Leites\\               
Department of Mathematics\\
Stockholm University\\
Albanov\"agen 28\\
SE-114 19 Stockholm\\
Sweden\\
dimleites@gmail.com
}

\begin{document}


\maketitle

\begin{abstract}
(A) The conventional integration theory on supermanifolds had been
constructed in order to have (an analog of) the Stokes formula in which a~sub-supermanifold is of codimension $1=(1|0)$. 
I review other integrations and formulate related open problems:  

(1) On the $1|1$-dimensional superstring associated with the trivial bundle,  in presence of a~contact structure there is a special integration useful in describing super versions of elliptic functions. It is needed to construct a~particular spinor representation of the Neveu-Schwarz superalgebra.

(2) Versions of the Stokes formula with ``over-supermanifold" of codimension $(0|-1)$ due to Shander and Palamodov should be developed further.

(3) Apply Shander's integration with odd parameters over chains to inverse problems.

(4) Establish existence of conjectural integrations (apparently, not leading to any analog of the Stokes
formula) related to
various (super)traces on various Lie superalgebras and the corresponding (super)determinants.

(B) I offer analogs of the cross
ratio for ``classical superspaces'', including infinite-dimensional versions. For four solutions of a~super matrix Riccati equation, these invariants are first integrals.
\end{abstract}

This review unites several ostensibly unrelated topics and their \lq\lq analytical continuations" and draws attention to open problems I find rather important: their solutions would lead us further to unpredictable results. I also wanted to give adequate references in certain cases under-represented in the literature and not as well-known as they deserve to be as, I hope, I demonstrated here why. The text contains definitions, construction, and problems (their solutions are future theorems).

For the basics on Linear Algebra in Superspaces and Supermanifold
theory, see Bernstein's lectures in the book \cite{D} (or \cite{L, B} if you read Russian); for some useful facts, see \cite{Sh14, Ser}. Of several sources listing the simple (and close to them) Lie superalgebras of finite dimension or $\Zee$-graded of polynomial growth, I recommend most recent \cite{GLS2, CCLL, BLLoS} with new results and insight. The list of these Lie superalgebras is rather long, in order not to reproduce it, I refer to \cite{BGLLS}.

\section{Integration}

\subsection{Integrations with the Stokes formula} In recent notes \cite{W}, Witten recalls integration suggested in the mid-1970's by
J.~Bernstein and me (see \cite{BL1, BL2}) and lists several applications. Here, I overview what is not covered in \cite{W}, and --- a bit differently --- what is. In \cite{BL1, BL2}, we offered an analog of
integration theory on supermanifolds superizing the Stokes formula. We had at our disposal the following two ingredients: 

(i1) the \textit{differential forms}, i.e., polynomials in the differentials $\hat x_i:=dx_i$
of the coordinates on the given supermanifold; the \textit{degree}-$k$ forms span the module $\Omega^k$ over the algebra $\cF$ of functions.   

(i2) \textit{volume forms}, the latter constituted a~rank
1 module $\Vol$ over $\cF$; 
\[
\begin{minipage}[l]{13.5 cm}
\textsl{under the change
of coordinates $x\mapsto y$ the generator $\vvol(y(x))$ of $\Vol$ is multiplied by the
Berezinian (superdeterminant) of the Jacobi matrix} $\frac{\partial y}{\partial x}$.
\end{minipage}
\]
Each of the above notions (i1) and (i2) had to be carefully reconsidered in super
setting because even the most innocent-looking notions and theorems
 displayed, in the super setting, inexplicable --- at first glance ---
signs at unexpected places (e.g., the Fubini theorem), see \cite{L}.


We wanted to
have some analog of the Stokes formula, so we needed two more ingredients: 

(i3) Certain integrands of ``degree'' smaller than that of volume forms and a differential connecting integrands of degrees that differ by 1. 

(i4) The notion of \textit{supermanifold-with-boundary}, see 
\cite{BL1}. (Observe that like manifold-with-boundary is not a~manifold, supermanifold-with-boundary is not a~supermanifold.)

Speaking of ``something'' to serve as integrand on a subsupermanifold we naturally wish to label these somethings by bidegree corresponding to the superdimension of the subsupermanifold.

To have an integration theory, one needs not only an object to be
integrated (the \textit{integrand}), but over what it is integrated (the \textit{cycle}), and the \textit{orientation}.
The latter two notions turned out to be more involved than we
originally thought; Shander clarified this in his development of
integration theory, see \cite{Sh2} expounded in \cite[Ch.~6]{L}. In particular, Shander introduced the objects dual to cycles and also labelled by bidegree (as superdimension).

Actually, having defined \textit{supermanifold-with-boundary}, we had enough to construct the
integration theory with an analog of the Stokes formula: by setting
$\Sigma_{-i}:=\Hom_{\cF}(\Omega^i, \Vol)$ we obtain a~complex dual to the de Rham complex
with $\Sigma_{0}= \Vol$ as forms of the highest degree. J.~
Bernstein and I called the elements of $\Sigma_{\bcdot}=\oplus \Sigma_{-i}$ \textit{
integraBle} forms (the ones one can integrate) and explained how to
integrate such forms in \cite{BL1}. More precisely, recall that $\Omega^i:=T(E^i(V^*))$, where $V$ is the tangent space, $E^i$ is the operator of raising to the $i$th exterior power, and $T(W)$ designates the space of $W$-valued tensor fields on the $(a|b)$-dimensional superdomain $\cU$ for every $\fgl(a|b)$-module $W$. Then, $\Sigma_{-i}:=T(E^i(V)\otimes \str)$, where $\str$ denotes the supertrace on $\fgl(a|b)$ and the module it defines; clearly,
$\Vol:=T(\str)$. Thus, locally, $\Sigma_{-i}$ is the $\Vol$-valued $\cF$-module of polynomials in $\partial_{\hat x_i}$; no generalized functions are needed to describe integrable forms, cf. \cite{W}.

(i5) The analog of exterior differential  $d:\Sigma_{-1}\tto\Vol$. To find such an operator and \textit{investigate if it is unique}, we had to consider the super analog of Veblen's problem solved, as on cue, by Rudakov shortly before, thus giving an invariant characterization of the exterior differential,  spaces on which it is defined, and pointing a~way to find its possible superanalogs.

\subsubsection{Veblen's problem and Rudakov's solution of it in the unary case} The exterior differential $d$ (instrumental in the Stokes formula)
is, as is well known, an invariant operator; in \cite{R} it is proven that $d$ is
the \textbf{only} invariant unary differential operator between spaces $T(V)$ of
tensor fields --- sections of vector bundles (over the $n$-dimensional manifold) whose fibers are irreducible
$\mathfrak{gl}(n)$-modules $V$ with a~vacuum vector. Remarkably, among all spaces
$T(V)$ 
\be\label{finDim}
\text{only the spaces 
$T(E^i(\id^*))=\Omega^i$ \textup{(}their fibers are finite-dimensional\textup{)}, }
\ee
where $\id$ is the identity (tautological) $\mathfrak{gl}(n)$-module, have a non-trivial submodule: $\Ker d$ for $i<n$ and $\mathrm{Im}\, d$ for $i=n$. 

In order to overview the options on supermanifolds, it was necessary to describe all differential operators $T(V)\tto T(W)$ invariant with respect to arbitrary (invertible) changes of coordinates --- obtain a~super version of Rudakov's result.

On the $n$-dimensional manifold (locally: a domain) $M$, let $W$ and $V_i$, where $i=1,\dots, k$, be irreducible $\fgl(n)$-modules. The problem 
\be\label{VP}
\begin{minipage}[l]{12 cm}
\textsl{Describe all $k$-linear differential operators} 
\[
D:T(V_1)\otimes_\cF T(V_2)\otimes_\cF \dots\otimes_\cF T(V_k)\tto T(W)
\]
\textsl{invariant with respect to any changes of
coordinates on  $M$}.
\end{minipage}
\ee
was raised in 1928 by O.~Veblen; for a~review, see \cite{GLS1} and an update \cite{BLS}. For unary ($k=1$ in problem \eqref{VP}) operators on manifolds, it was solved by Rudakov \cite{R} as a~part
of another problem, namely, the description of continuous in the discrete (or, dual to it, $x$-adic) topology irreducible vacuum
vector modules over the simple Lie algebra $\fvect(n|0)$ of formal or polynomial
vector fields.

\textbf{Aside}. An interesting and curious fact: the solution of the binary ($k=2$ in problem \eqref{VP}) version of Veblen's problem on manifolds (no ``super''), i.e., classification of invariant bilinear differential operators, due to Grozman \cite{Gr}, brings several analogs of the Poisson bracket on Lie \textbf{super}algebras close to simple. (Up to dualization and permutation of arguments, there are no other invariant bilinear operators.)

\subsection{Integrations unconventional: either without Stokes formula, or for another reason} For a graphic and short description of $\cC$-points in the non-super setting, see \cite{MaAG}. 

\sssec{Preliminaries: Automorphisms that preserve parity and their odd parameters} All parity preserving automorphisms over the ground field $\Kee$ of the Grassmann algebra $\Lambda^{\bcdot} (n)$  on $n$ generators $\xi_i$ are of
the form (see \cite{Dj})
\begin{equation}
\label{7eq8}
\xi_j\;
\mapsto\;\sum_{r} \varphi_j^{r}\xi_{r}+ \sum_{s\geq 1}\;
\sum_{j_1<\dots <j_{2s+1}} \varphi_j^{j_1\dots
j_{2s+1}}\xi_{j_{1}}\dots\xi_{j_{2s+1}},\text{~~where $\varphi_j^{r}, \varphi_j^{j_1\dots
j_{2s+1}}\in\Kee$,}
\end{equation}
with invertible matrix $(\varphi_j^{r})$. Clearly, such
automorphisms constitute the algebraic group (or a~Lie group if $\Char\Kee=0$) whose Lie algebra consists of the even
elements of the Lie superalgebra $\fvect(0|n):=\fder\, \Lambda^{\bcdot} (n)$ of superderivations of $ \Lambda^{\bcdot} (n)$.

Let $\cU:=(U, C^\infty(\cU):=C^\infty(U)\otimes \Lambda^{\bcdot} (n))$, where $U\subset \Ree^m$ is a domain, be an $m|n$-dimensional \textit{superdomain} --- the local chart of a~ supermanifold. \textbf{To account for the odd parameters (the boxed terms in formula \eqref{7eq99}) we have to consider the functor $\cC\longmapsto \Aut_\cC^{\ev}(C^\infty(\cU)\otimes \cC)$ from the category of supercommutative superalgebras to the category of groups}, sending every $\cC$  to the group of the parity preserving (as indicated by the superscript $\ev$) $\cC$-linear automorphisms of the form (here the $u_i$ are the coordinates in $U$, and $u:=(u_1, \dots, u_m)$)
\begin{equation}
\label{7eq99}
 \left\{\renewcommand{\arraystretch}{1.6}
 \begin{array}{ll}
 \varphi^*(u_i) \;& = \; \varphi^0_i(u) + \sum\limits_{r\geq 1}\;
 \sum\limits_{i_1<\dots<i_{2r}}\varphi_i^{i_1\dots
 i_{2r}}(u)\xi_{i_{1}}\dots \xi_{i_{2r}}\\
 &+
 \fbox{$\sum\limits_{r\geq 0}\;
 \sum\limits_{i_1<\dots<i_{2r+1}}\varphi_i^{i_1\dots
 i_{2r+1}}(u)\xi_{i_{1}}\dots \xi_{i_{2r+1}}$}, \\
\varphi^*(\xi_j) \;& = \sum\limits_{r\geq 0}\; \sum\limits_{j_1<\dots <j_{2r+1}}
 \psi_j^{j_1\dots
 j_{2r+1}}(u)\xi_{j_{1}}\dots\xi_{j_{2r+1}}\\
&+\fbox{$\psi^0_j(u)+
\sum\limits_{r\geq 1}\; \sum\limits_{j_1<\dots <j_{2r}}
 \psi_j^{j_1\dots
 j_{2r}}(u)\xi_{j_{1}}\dots\xi_{j_{2r}}$},\\
\end{array}
\right.
\end{equation}
where for all $r$, the elements
\[
\text{$\varphi^0_i(u),\varphi_i^{i_1\dots
 i_{2r}}(u), \ \psi_j^{j_1\dots
 j_{2r+1}}(u)\in \cC_\ev$ and 
$\psi^0_i(u),\psi_j^{j_1\dots
 j_{2r}}(u),\ \varphi_i^{i_1\dots
 i_{2r+1}}(u)\in \cC_\od$ }
 \]
are the parameters of the infinite-dimensional supergroup of automorphisms of $C^\infty(\cU)$ or, equivalently, of diffeomorphisms of the superdomain $\cU$; infinitesimally, these are parameters of the Lie superalgebra $\fvect(m|n)$. The automorphisms \eqref{7eq99} are called \textit{$\cC$-points} of the supergroup of diffeomorphisms of $\cU$; the ones without the boxed terms, i.e., the $\cC$-points for $\cC=\Kee$, constitute the underlying group of this supergroup.

\sssec{Automorphisms that do not preserve parity} In 1973 or so, I've astonished Berezin informing him that there are automorphisms not preserving parity of the Grassmann algebra $\Lambda^{\bcdot} (n):=\Cee[\xi]$, where $\xi=(\xi_1,\dots, \xi_n)$, considered as an associative algebra, cf. \cite{Be}, where all automorphisms are supposed to preserve parity. The first complete description of all automorphisms of $\Lambda^{\bcdot} (n)$ not preserving parity of the Grassmann algebra over a field of any characteristic is due to Djokovi\v{c} (see \cite{Dj}).

The group $\Aut (\Lambda^{\bcdot} (n))$ of all automorphisms of $\Lambda^{\bcdot} (n)$ is just a bit larger than the group $\Aut^{\ev}(\Lambda^{\bcdot} (n))$ of parity preserving automorphisms: the general automorphism $\varphi$ of $\Lambda^{\bcdot} (n)$ is of the form
\be\label{mix}
\varphi(\xi_i)=\varphi_0(\xi_i)(1+\psi(\xi))\text{~~ for all $i$ and some $\psi(\xi)\in\Lambda^{\bcdot} (n)_\od$},
\ee
where $\varphi_0$ is a parity preserving automorphism.

More generally, for the classification of non-preserving parity automorphisms of $\Cee[x,\xi]$, where $x=(x_1,\dots, x_m)$, see \cite{LS1, LS2} and for details, see \cite{Ba}. The meaning of such general automorphisms is unclear at the moment, see \cite{Lq}. 

\sssec{Towards conjectural interpretations} By analogy with supersymmetries --- which are wider than symmetries carried by groups or Lie algebras --- these general automorphisms that do not preserve parity widen supersymmetries still further. More precisely, the even functions are invariantly defined (bosons do not depend on the coordinate systems), whereas to say that a~ function is odd (corresponding to a~fermion) is the same as to select a~coordinate system modulo the subspace of even functions. \textbf{The notion of Fermion is not invariant}.

Recently, U.~Iyer proved that the non-homogeneous with respect to parity subalgebras of Lie superalgebras, called \textit{Volichenko algebras}, play the role of Lie algebras for the groups of $\cC$-points of such general automorphisms of supercommutative superalgebras, see \cite{I}.

\sssbegin{Theorem} For $n=2k+1$, the Berezin integral is invariant with respect to the general changes of coordinates \eqref{mix}: the Berezin integrals of $\varphi(f)$ and $\varphi_0(f)$ coincide for every $f\in\Lambda^{\bcdot} (n)$. For $n=2k$, this holds for every even $f$, but not for the odd ones. \end{Theorem} 

\begin{proof} The integral is equal to the coefficient of the highest term of the integrand. Set $a_i:=\varphi_0(\xi_i)$. Since $\psi$ is odd, $(1+\psi)a_j=a_j(1-\psi)$, and since $\psi^2=0$, we have
\[
\varphi(\xi_i)\varphi(\xi_j)=a_i(1+\psi)a_j(1+\psi)=a_ia_j(1-\psi^2)=a_ia_j.
\]
Hence, $\varphi(f)=\varphi_0(f)$ for $f$ even and $\varphi(f)=\varphi_0(f)(1+\psi)$ for $f$ odd. For $f$ odd, every homogeneous component of $\varphi_0(f)\psi$ is of even degree, so it contributes nothing to the highest term if $n$ is odd. If $n=2$, then $\varphi_0=\id$ and $\psi=\xi_2$ give $\varphi(\xi_1)=\xi_1+\xi_1\xi_2$, whose integral is $1$, whereas that of $\xi_1$ is $0$.
\end{proof}

\ssec{Miscellanea}{}~{}

\sssec{Non-trivial densities on supermanifolds with an odd non-degenerate differential 2-form (periplectic form)} For these, see H.~Khudaverdian's papers \cite{H, H2} and \cite[Ch.4, \S7]{Ma}.

\sssec{No invariant differential operator depends on higher jets of diffeomorphisms} Ruda\-kov proved that there are no (non-scalar) invariant differential operators between spaces of objects depending on higher jets of the diffeomorphism (changes of coordinates); the proof of this fact constitutes the bulk of \cite{R}. 

A.~Schwarz (a.k.a. Shwarts and \v{S}warc) et al. 
attempted to integrate
densities and objects depending on higher jets of the diffeomorphism,
but all their non-trivial examples boil down to either pseudo-differential
(\cite{BL2}) or integrable forms (depending only on the 1-jets of diffeomorphisms). These \textit{ad hoc} attempts suggest the following

\sssec{Problem} Prove that on supermanifolds, same as on manifolds, there are no invariant differential operators depending on higher jets of diffeomorphisms. 
Observe that no \textit{direct} superization of Rudakov's proof of the above problem on manifolds is visible; it seems, one needs totally new ideas.

\sssec{Analogs of the Berezin integral of tensors constituting irreducible infinite-dimen\-sional fibers} Having obtained an analog of
Rudakov's result for the general vectorial Lie \textit{super}algebra $\fvect(a|b)$, see
\cite{BL3}, we can be sure that 
\be\label{(*)}
\begin{minipage}[l]{12 cm}
there is only one integration theory
on supermanifolds\\
 \textbf{provided the integration involves only
tensors\\ with irreducible finite-dimensional spaces}.
\end{minipage}
\ee

The result of \cite{BL3} does not preclude, however, unconventional
integrations with infinite-dimensional fibers. In
\cite{GLS1} the calculations from \cite{BL3} are used to classify invariant differential operators between sections with
infinite-dimensional fibers snubbed at in \cite{BL3} for no reason
except tradition: in the non-super case, the invariant operator exists only if fibers of tensor fields are finite-dimensional, see eq.~\eqref{finDim}. 

Having stated that it is impossible to integrate
differential forms on supermanifolds (functions polynomial in differentials of coordinates) in \cite{BL1}, Bernstein and I showed in \cite{BL2} how to integrate \textbf{pseudo}-differential forms, i.e., functions \textbf{smooth} in
differentials of coordinates: $\hat x:=dx$. Observe that on supermanifolds of superdimension $(a|b)$, the space of pseudo-differential forms of a given degree $\lambda\in\Kee$ is a module of finite rank over the algebra $\cF$ of functions only for $b=0$ or 1.

In \cite[Ch.6]{L}, Shander singled out several interesting types of pseudo-differential forms.
For example, one can consider \textit{homogeneous} forms, i.e.,
such that 
\[
f(x, t\hat x)= t^a f(x, \hat x) \text{~~for all $t\in\Ree$ and some $a\in\Ree$};
\]
or \textit{forms with compact support in} $\hat x_i$ for every odd indeterminate $x_i$;
or \textit{forms rapidly decreasing as} $\hat x_i\tto \infty$ for all odd $x_i$, etc.

\sssec{Problem} Compute 
\textit{semi-infinite cohomology of supermanifolds} related to semi-infinite pseudo-differential forms.

\textup{No examples of absolutely new invariants like this semi-infinite cohomology are calculated yet, except for examples due to A.~Zorich (see \cite{Z1}, and a~general discussion in \cite{VZ}). For Shander's clarification of bigraded differential forms introduced in \cite{VZ} and \textbf{how to integrate them with the help of odd parameters}, see \cite[Ch.~6]{L}.}

\textup{Observe that several authors published some computations of bigraded cohomology on supermanifolds, but these papers contain irreparable mistakes as explained, e.g.,  in \cite{Kato}. }

\ssec{A version of the Stokes formula with codimension $(0|-1)$ ``over-supermanifold"}
Around 1986, Shander wrote ca 80 pages elucidating $2+2$ pages of \cite{BL1,BL2}. Shander's results were published ca 25 years later and so far only in Russian, see \cite[Ch.~6]{L}. There, Shander demonstrated the importance of odd parameters (see also \cite{L2}) and implicitly suggested a~version of the Stokes formula with ``oversupermanifold" of codimension $(0|-1)$. 

In 1996, Palamodov published a~note \cite{Pal}, where he gave an example of the Stokes formula in which the ``boundary'' of the $0|0$-dimensional point is a~$0|1$-dimensional line.

\sssec{Problem} Develop integration via the ideas due to Shander and Palamodov. 

\ssec{Supertraces and superdeterminants} Since 1975,
I wondered what if we stop insisting on having an analog of the
Stokes formula? What remains of the paraphernalia needed for  integration then? 

The superdeterminant of the Jacobian, i.e., of the matrix $\left(\frac{\partial y}{\partial x}\right)$ of the change of local coordinates on the supermanifold $\cM$ will definitely remain. Since it is easier to deal with Lie
superalgebras than with supergroups, let me list analogs of trace for Lie
superalgebras. Then, the Lie superalgebra $\mathfrak{g}$ can be
exponentiated to a~Lie supergroup, and we can consider the super determinant defined (on the level of points over a~supercommutative superalgebra $\cC$) via the formula
\be\label{(1)}
\det \exp X= e^{\tr (X)} \quad\text{ for every } \; X\in(\mathfrak{g}\otimes \cC)_\ev. 
\ee
The same procedure can be applied to any $\Zee$-graded Lie superalgebra $\fg$ of vector fields on whose component $\fg_0$ a~ trace is defined. Let us consider examples.

A) Let $\cM$ be a~supermanifold with the fiber of the tangent bundle isomorphic to the space $V$; let $\mathfrak{g}$ be a~Lie subsuperalgebra of $\fgl(V)$ with its own trace
$\tr$ (which by definition is \textit{any linear function on $\mathfrak{g}$ such that $\tr([x, y])=0$ for all
$x,y\in\mathfrak{g}$}; in the super setting, this trace is often called \textit{super}trace for emphasis). Then, for the role of $\Vol$ over $\cM$, we can take tensor
fields of type $\tr$, i.e., take sections of rank 1 module given by $\tr$ over the algebra $\cF$ of functions on $\cM$; the change of variables acts on $\Vol$ via the Lie derivative. The Lie superalgebra of infinitesimal transformations of this $\Vol$ is the
\textit{Cartan prolongation} (for its definitions, see \cite{Shch, GLS1, BGLLS}) of the pair $(V,
\mathfrak{g})$, where $V$ --- the fiber of the tangent bundle on $\cM$ --- is the tautological
module over $\mathfrak{g}\subset\fgl(V)$. 

The prime example of the Lie superalgebra with trace analogous to $\fgl(V)$ is provided by the Poisson Lie superalgebra
$\mathfrak{g}=\mathfrak{po}(0|2n)$. Indeed, there is a~parametric family
(deformation, physicists called it \textit{quantization}) of Lie superalgebras $\mathfrak{g}_{t}$ which at $t=0$
coincides with $\mathfrak{po}(0|2n)$ and $\mathfrak{
g}_{t}\simeq\mathfrak{gl}(2^{n-1}|2^{n-1})$ for every $t\neq 0$, see
\cite{LSh2}. Under quantization the trace given on $\mathfrak{po}(0|2n)$ by the Berezin integral turns into the supertrace on $\mathfrak{gl}(2^{n-1}|2^{n-1})$.

B) From various points of view $\mathfrak{gl}(n)$ has
at least two superanalogs: the ``simple-minded'' superanalog,
$\mathfrak{gl}(n|m)$, and the ``queer'' superanalog, $\mathfrak{q}(n)$ which corresponds to the case of the Schur lemma with irreducible representations that have no invariant sub\textbf{super}spaces, but have an invariant subspace. It is easy to show that over an algebraically closed field $\Kee$ such a~representation commutes with an odd operator $J$ such that $J^2=\lambda \id$, where $\lambda \in\Kee^{\times}$ and $\sdim V=(n|n)$, the invariant subspace being an eigenspace of $J$. If the characteristic of $\Kee$ is $\neq 2$, it is possible to select a~basis in $V$ so that $J=\begin{pmatrix}0&1_n\\-1_n&0\end{pmatrix}$, and then we can realize the supermatrices of $\fq(n):=\{X\in \fgl(n|n)\mid [X, J]=0\}$ in the form $(A,B):=\begin{pmatrix}A&B\\ B&A\end{pmatrix}$, where $A,B\in\fgl(n)$. On
$\mathfrak{q}(n)$, the supertrace vanishes identically, but there is
the \textit{queer trace} $\qtr: (A,B)\mapsto \tr B$. On the group $GQ(n; \cC)$ of invertible matrices of the form $\begin{pmatrix}A&B\\ B&A\end{pmatrix}$, where $A\in\fgl(n;\cC_\ev)$ and $B \in\fgl(n;\cC_\od)$, where $\cC$ is a supercommutative superalgebra,  the \textit{queer determinant} is defined by the formula
\[
\qet(A, B):=1+\tau\sum_{i\geq 1}\frac{\tr((A^{-1}B)^{2i-1})}{2i-1}\in\cC[\tau], \text{where $p(\tau)=\od$}.
\]
The queer trace is odd, and therefore \textbf{to
describe the representation corresponding to it, and the queer determinant $\qet$, we need an odd parameter} $\tau$, see \cite[Ch.~1]{L}; compare with the first definitions in \cite{BL4}.

C) Still other traces are
related to the ``quasiclassical limit'' of $\str$ and $\qtr$ as
$t\longrightarrow 0$: the restriction of the above quantization is a
parametric family $\mathfrak{g}_{t}$ which at $t=0$ coincides with
$\mathfrak{po}(0|2n-1)$ and $\mathfrak{g}_{t}\simeq\mathfrak{
q}(2^{n-1})$ for $t\neq 0$, see \cite{LSh2}. 

Under quantization the trace given on $\mathfrak{po}(0|n)$ by the Berezin integral turns into $\qtr$ on $\mathfrak{q}(2^{k-1})$ for $n=2k-1$, and into $\str$ on $\mathfrak{gl}(2^{k-1}|2^{k-1})$ for $n=2k$.

\ssec{Integration on $1|1$-dimensional superstrings with contact structure} Consider the form $\alpha:= dt + \tau d\tau$ on real or complex supermanifold $\cM=(M,\cO)$ of dimension $1|1$ with coordinates $t$ (even) and $\tau$ (odd). Let $\fk^L(1|1)$ be the Lie superalgebra of contact vector fields preserving the distribution singled out by $\alpha$; the superscript $L$ is for ``Laurent" since $\fk^L(1|1)$ can be considered as a subsuperalgebra in the Lie superalgebra of vector fields with Laurent coefficients (in $t$). All tensor bundles on the spaces of whose sections the Lie superalgebra $\fk^L(1|1)$ acts are glued of the simplest ones --- linear (of rank $1=1|0$ or $0|1$) which are labeled by their ``conformal weight'' --- the eigenvalue of the generator under the action of the vector field $2t\partial_t+\tau\partial_\tau$. The space of differential 1-forms does not have $\fk^L(1|1)$-invariant basis over the algebra of functions $\cF$, only the subspace $\cF\alpha$ is  invariant; the quotient modulo $\cF\alpha$ is isomorphic to the space of volume forms, i.e., $[d\tau]=\vvol(t,\tau)$. Hence, we can integrate not only functions, but also 1-forms (with compact support):
\[
\begin{array}{l}
\int_{\cM} (fdt+gd\tau)=\int_{\cM} (f(dt+\tau d\tau -\tau d\tau)+gd\tau):=\\
\int_{\cM} (g-f\tau)[d\tau]=\int_{\cM} (g-f\tau)\vvol(t,\tau)=\int_{M}(g_1-f_0) dt, 
\end{array}
\]
for all functions $f=f_0(t)+f_1(t)\tau$ and $g=g_0(t)+g_1(t)\tau$. This construction, although is based on Berezin's integration, is another integration, existing in parallel with Berezin's. It is used in Levin's theory of elliptic functions on supermanifolds; for an exposition (in which the above-described mechanism is not explained since it was only a tool, not the object), see \cite{Mn}.

\sssec{Problem} Construct explicitly \textup{(by giving formulas for the operators)} the spinor representation \textup{(see \cite{LSh1})} of $\fk^L(1|1)$ related to the anti-symmetric non-degenerate bilinear form on $(\cF\otimes \cC)/\cC\cdot 1$ \textup{(cf. \cite{LSa})} given by
\[
\begin{array}{l}
(f, g):=\int_{\cM} fdg:=\int_{\cM} f(dt\partial_tg+d\tau\partial_\tau g)\\
=\int_{\cM} ((-1)^{p(f)}(dt +\tau d\tau -\tau d\tau)f\partial_tg+d\tau f\partial_\tau g)\\
=\int_{\cM}(f\partial_\tau g -(-1)^{p(f)}\tau f\partial_tg)\vvol(t,\tau)=\int_{M}(f_1g_1-(-1)^{p(g)}f_0\partial_tg_0) dt.
\end{array}
\]

\ssec{TKK construction} In the 1960s, a~remarkable correspondence between Jordan algebras and
certain $\Zee$-graded Lie algebras became explicit, cf. \cite{T},
\cite{K1}--\cite{K2} and \cite{Koe}. Kantor used this
correspondence to list simple Jordan algebras (over $\Cee$ and
$\Ree$) by the, so far, simplest known method. He clarified the
mysterious relation of Jordan algebras with classical domains and
actively studied certain generalizations of Jordan algebras
associated with $\Zee$-graded simple Lie algebras $\mathfrak{g}=\mathop{\oplus}\limits_{-d\leq i\leq h}\mathfrak{g}_i$ of finite depth
$d$. (Since all finite-dimensional simple Lie algebras are ``symmetric'',
$\mathfrak{g}_i=(\mathfrak{g}_{-i})^*$,
their \textit{height} $h$ is equal to $d$.) The Lie superalgebras of supergravity and supertwistors also have symmetric $\Zee$-gradings but of depth $d>1$, see \cite{Ma, GL}; so it is interesting to find generalizations of Jordan
algebras (or, rather, any useful structures) related with such gradings.

Following Freudenthal and Springer, Kantor generalized products
\eqref{(2)} to several arguments; this is natural for $d>1$. 

I suggest,
contrariwise, to stick to formula \eqref{(2)}, even for $d>1$, with the
minimal modification: fix $j\in\mathfrak{g}_{1}$ and for $\mathfrak{
g}_{-}=\mathop{\oplus}\limits_{i< 0}\mathfrak{g}_i$ set 
\be\label{(3)}
x\circ y= [[j, x], y]\quad\hbox{ for all } x, y\in
\mathfrak{g}_{-}. 
\ee
In this way, we obtain non-commutative generalizations of Jordan
algebras. Superization is immediate: via the Sign Rule. 

\sssec{Problem} What are the relations additional to the Jordan identity \eqref{(JoI)} (and commutativity or replacing it) these generalizations satisfy? What type of integrable systems are
associated with them in the Sokolov-Svinolupov approach, cf.
\cite{SS,HSY}? 

Tables in \S\ref{Tab} (borrowed from \cite{LSV} and completed with an exceptional example discovered by Shchepochkina) provide us with a~list of simple
Jordan superalgebras associated with the known 
simple
$\Zee$-graded Lie superalgebras of polynomial growth (SZGLSAPGs for
short), cf. \cite{Ka}, \cite{LSh2}, \cite{KMZ, CaKa}, namely, with
$\Zee$-gradings of depth 1 of SZGLSAPGs \textit{including} finite-dimensional ones. 

\subsubsection{The Tits--Kantor--K\"ocher functor $\mathfrak{kan}$ }

Let $J$ be a~Jordan algebra and $j$ the tensor that determines
the product in $J$, i.e., $j(x, y) := x \circ y$. To $J$, Kantor
assigned (see \cite{Ka} with its epigraph ``There are no Jordan algebras, there are only Lie algebras. I.~Kantor") a~$\Zee$-graded Lie algebra
$\mathfrak{kan}(J):=\mathop{\oplus}\limits_{|i|\leq 1}
\mathfrak{kan}(J)_{i}$, a~$\Zee$-graded Lie subalgebra in the Lie algebra 
$\mathfrak{vect} (J)$ of vector fields on $J$ in the standard $\Zee$-grading ($\deg x_i=1$) such that 
\[
\renewcommand{\arraystretch}{1.4}
\begin{array}{l}
\mathfrak{kan}(J)_{-1}:=\mathfrak{vect} (J)_{-1}, \\
\mathfrak{kan}(J)_{0}:= \Span(L_{a}, [L_a, L_b]\mid a, b \in J), \text{~~where $L_a(x):= a\circ x$},\\
\mathfrak{kan}(J)_{1}:= \Span(j, [L_a, j]\mid \;a\in J).
\end{array}
\]

There are versions of this construction due to Tits and K\"ocher which yield the same Lie algebra when the constructions are applied to a~simple finite-dimensional Jordan algebra. But, when these constructions are applied to Jordan \textbf{super}algebras, they can produce non-isomorphic Lie superalgebras, as demonstrated in \cite{BC}. However, Kantor's idea was that Jordan algebras are determined by Lie algebras, so this result in \cite{BC} is opposite Kantor's preferred direction of study. 

\subsubsection{Digression: on homological fields}
I do not know what structure is related with an arbitrary odd
element $L\in\fg_1$ of the $\Zee$-graded Lie superalgebra $\fg = \oplus_{i\geq -1}\ \fg_i$, but if $L$ is \textit{homological}, i.e., $[L,L]=2L^2=0$,
then the formula
\be\label{(L)}
[x, y]' :=[[L, x], y] 
\ee
determines a~Lie superalgebra structure on $\Pi(\mathfrak{g}_{-1})$.
This structure had been first noticed, perhaps, by M.~Gerstenhaber
in the 1960s and rediscovered many times since then. It is
interesting to describe in intrinsic terms the $L$'s which determine
\textit{simple} Lie (super)algebras.

\textit{Non-degenerate} homological  vector fields were first characterized, in connection with
the problem of integration of differential equations on
supermanifolds, by V.~Shander, see \cite{Sh1}, who gave their~normal forms. 

Shander did not consider fields with
singularities; 
Vaintrob did: see, e.g., \cite{V}. Vaintrob and Serganova showed
that the study of singularities of homological fields, and their
classification, turns out to be rather similar to the case of
singularities of smooth functions, see \cite{SW}. Regrettably, the answer describing the elements $L$
corresponding to the simple Lie (super)algebras is more trivial than expected,
as I have established with Grozman:

Let, first, $\mathfrak{g}=\mathfrak{gl}(n)$; denote the matrix units
by $\partial_{i}^j$; let $(\partial_{i}^j)^*=x_{j}^i$ be the dual
basis. Then, $L=\sum \xi_{j}^i\xi_{k}^j\frac{\partial}{\partial\xi_{k}^i}$, where
$\xi_{j}^i$ is the odd copy of $x_{j}^i$, is homological. Similarly, if the $X_{i}$ form
a basis of $\mathfrak{g}$ and $[X_{i}, X_{j}]=\sum c_{ij}^kX_{k}$,
then the operator $L\in \mathfrak{vect}(0|\dim \mathfrak{g})$ of
the form 
$\sum c_{ij}^kX_{i}^*X_{j}^* \partial_{X_{k}^*}$, where $X_{i}^*$
is the \textup{(}odd\textup{)} dual of $X_{i}$, is homological. 

Having observed that every simple finite-dimensional Lie algebra
(over $\Cee$) has a~non-degenerate invariant symmetric bilinear form (NIS), we
see that, for any such Lie algebra $\fg$, the vector field $L$ is hamiltonian;
 its generating function is, up to a factor, the Cartan $3$-form $(x,y,z)\mapsto B([x,y],z)$, where $B$ is the NIS; it is of degree $1$ in the standard grading of the Poisson Lie superalgebra $\mathfrak{po}(0|\dim \mathfrak{g})$ and $\fg$-invariant, in particular of weight $0$ with respect to the maximal toral subalgebra of $\fg\subset\fo(\dim\fg)=\mathfrak{po}(0|\dim \mathfrak{g})_0$. 

Generalization to Lie superalgebras is straightforward. Still,
observe that some simple Lie superalgebras have no NIS at all, some
(e.g., $\fpsq(n)$) have an \textit{odd} NIS, in which case $L$ belongs to the
Buttin superalgebra (with antibracket a.k.a. the Schouten bracket), not to the Poisson superalgebra. (For the classification of NISes on finite-dimensional Lie (super)algebras, and infinite-dimensional Lie (super)algebras of polynomial growth, see \cite{BKLS}.)

\subsection{SZGLSAPGs of depth 1 and height 1}
All possible $\Zee$-gradings of\linebreak
 SZGLSAPGs $\mathfrak{g}$ are listed in
\cite{Ka} for $\dim \mathfrak{g} < \infty$ and in \cite{L0} for 
the other cases known in 1982. Our job is to pick those of them which are of depth
1, in particular, of the form $\mathop{\oplus}\limits_{|i| \leq 1}\;
\mathfrak{g}_i$, see Tables in \S\ref{Tab}.

Observe that Albert's notation for Jordan algebras follows
Cartan's notation for the corresponding Lie algebras. As follows
from Serganova's classification of systems of simple roots of simple
Lie superalgebras \cite{S}, Cartan's notations are highly
inappropriate for Lie superalgebras (hence, for Jordan superalgebras); for more arguments, see \cite{CCLL}.

\section{Cross ratios}

Here, I consider certain algebraic structures associated with several
selected ``classical superdomains''. So far, nobody knows yet (as far
as I know) what should be considered as such a~superdomain, except for the case of finite-dimensional manifolds when all is clear. For
supermanifolds, we just consider the easiest to handle simple Lie
supergroups and their Lie superalgebras whereas the elusive ``right''
definition requires, perhaps, ``semisimple" or ``almost simple" Lie
superalgebras (for their definition, see \cite{BGLLS}). So we take the road of least resistance:

The ground field is $\Cee$, the classical
superdomains are considered as quotients of \textit{simple} or close to
them ``quasi-reductive'' Lie supergroups (see \cite{Sqr}) modulo \textit{certain} maximal
parabolic subsupergroups; for their list, see Tables.
Consider the following questions:
\begin{equation}\label{Q}
\begin{minipage}[c]{13cm}
 (Q1) What are the criteria for selecting \textit{certain} sub(super)groups
among\\ 
\phantom{.}\hskip 1cm other maximal ones in the ambient supergroup?\\
 (Q2) Why cosets modulo other parabolic subsupergroups are seldom\\
\phantom{.}\hskip 1cm considered in Differential Geometry? Observe that precisely these\\ \phantom{.}\hskip 1cm ``other''
cosets are the main topic of study in analytical mechanics\\
\phantom{.}\hskip 1cm of non-holonomic dynamical systems and in supergravity, see\\ 
\phantom{.}\hskip 1cm \cite[Ch.5]{Ma}, \cite{GL, BGLS}.
\end{minipage}
\end{equation}

In these questions, ``super'' is beside the point, so we can very well
begin with manifolds. The classical domains are distinguished among
symmetric spaces by the fact that the Lie algebra of the symmetry
group of any classical domain (Hermitian symmetric space) $M$ is a
simple complex Lie algebra of the form
\be\label{(d=1)}
\mathfrak{g} = \mathop{\oplus}\limits_{|i|\leq 1}\;
\mathfrak{g}_i; 
\ee
the tangent space to $M$ at a~fixed point can be identified with the
$\mathfrak{g}_{-1}$ and, on it, one can always define a~Jordan
algebra structure, by fixing any element $j\in\mathfrak{g}_{1}$ and
setting
\be\label{(2)}
x\circ y:= [[j, x], y]\quad\hbox{ for all } x, y\in
\mathfrak{g}_{-1}. 
\ee
Jordan multiplication on the prime example --- the space of $n\times n$-matrices is no less natural than the Lie multiplication, 
but the identity \eqref{(JoI)} it satisfies is less known than the Jacobi identity. (No wonder: it has no graphic interpretation, unlike Jacobi identity which is just the Leibniz rule for the adjoint action $\ad$ or the bracket.) Recall that a~\textit{Jordan algebra} is a~commutative
algebra $J$ with product $\circ$ satisfying, instead of
associativity, the identity
\be\label{(JoI)}
(x^2\circ y)\circ x = x^2\circ (y\circ x). 
\ee

In the very inspiring paper and
book \cite{Mc, Mcb}, McCrimmon gave an account of some applications of
Jordan algebras from antiquity to nowadays. All numerous literature known to me
(\cite{Koe, Mc, U} and references therein) 
strengthens my prejudice
that \textit{general} Jordan algebras are not very useful. 
Contrariwise, {\bf simple} Jordan algebras are related to several
notions important in various problems. For simple Jordan algebras,
their \textit{generic norm} (see \cite{K3}) should be non-degenerate
and this imposes additional constraints on the parabolic subalgebra.
This answers the above questions \eqref{Q}, but tempts one to make use of the
other coset spaces spoken about in (Q2), as well.

I wish to make similar use of simple Jordan \textit{super}algebras,
especially infinite-dimensional ones, associated with infinite-dimensional superdomains listed in \cite{LSV} and Tables \S\ref{Tab}. 

In \cite{K3}, Kantor generalized the cross ratio of four points from a function on
$\Pee^1$ to a function on certain quotients $G/P$, where $G$ is a~simple Lie
group and $P$ is its parabolic subgroup. The Lie algebra
$\mathfrak{g}:=\hbox{Lie}(G)$ in these cases is of the form
$\mathfrak{g}=\mathop{\oplus}_{|i| \leq d}\; \mathfrak{g}_i$.

Kantor's construction of cross-ratio naturally
appears in a~ super setting \cite{Mn}; this prompts me to try to
decipher a~part of \cite{K3}; see also \cite{LiSS, San}.

I will consider here the simplest case, the case where $d=1$, and the
corresponding Jordan (super)algebra is simple. In this case, one can
generalize the cross ratio from $\Pee^1={\tt Gr}_{1}^2$ to a
collection of $\GL(2m|2n)$-invariants of four points on ${\tt
Gr}_{m|n}^{2m|2n}$. Kantor defined the cross ratio as
\be\label{(CR)}
(A, B, C, D):=(A-B)(C-B)^{-1}(C-D)(A-D)^{-1} \text{~~for all
$A, B, C, D\in J$}. 
\ee

If $mn=0$, then for any
$A, B, C, D\in (J\otimes \cC)_{\bar 0}$ and any supercommutative superalgebra $\cC$, I suggest to take for the \textit{cross
ratio} the collection of \textbf{all} coefficients of
\textbf{all} powers of $\lambda$: 
\be\label{CRd}
\det((A, B, C, D)-\lambda E), \text{~~where $E$ is the unit matrix}.
\ee

By dimension considerations these are all the invariants of four
points for the general, orthogonal and Lagrangian Grassmannians, but not for the quadrics.

If $mn\neq 0$, take the Berezinian instead of $\det$ in eq.~\eqref{CRd}. The amount of \textit{polynomially independent} invariants (coefficients of
the powers of $\lambda$) is
infinite, cf. \cite{Ser}, but if we consider rational
dependence, which is natural in super setting, the coefficients of
the first $n+m$ powers of $\lambda$ generate the algebra of
invariants and this collection is a~natural candidate for the cross ratio, or rather, the collection of all zeros and poles of $\ber((A, B, C, D)-\lambda E)$.

On $QGr_{n}^{m}$ which Manin called $\Pi$-symmetric Grassmannian $\Pi{\tt
Gr}_{n|n}^{m|m}$, see \cite[Ch.5,\S6, no.4]{Ma}, we should take the queerdeterminant, $\qet$,
instead of $\det$  in eq.~\eqref{CRd}; the collection of invariants thus obtained is finite, see \cite{Sh4}.

For (super)Grassmannians associated with loop (super)algebras and curved superGrassmannians, we have to consider
matrix-valued functions. In these two types of cases, formula~\eqref{CRd} returns a~collection of functions, rather than numbers.

\ssec{Related results} Four-point
functions in $N=2$ superconformal field theories were studied in \cite{EHPSW}. In \cite{Ze}, the
matrix cross ratio is applied to Riccati equation. In \cite{Ze1}, classical domains related to loop algebras are considered. See also a~beautiful paper \cite{MD}, where the super cross ratio of four points of the supercircle is constructed. Together with
\cite{Mn} these works provide a~wide range for applications of the new cross
ratios introduced above in the present article.

The invariants of points of the supercircle under $\OSp(1|2)$ became coordinates on the \textit{decorated super-Teichm\"uller space} of a~punctured surface: the even ones (super $\lambda$-lengths) and the odd ones (one for each triangle of an ideal triangulation), see \cite{PZ}; the transformation rules of these coordinates under a~flip of the triangulation are the super Ptolemy relations, see \cite{MOZ}; for $\OSp(2|2)$, see \cite{IPZ}.

\sssbegin{Proposition}[First integrals of the super matrix Riccati equation]\label{Ric} Let $\cC$ be a~supercommutative superalgebra, and let
\[
H(t)=\begin{pmatrix}\alpha(t)&\beta(t)\\ \gamma(t)&\delta(t)\end{pmatrix}
\]
be an even supermatrix of format $(m|n)+(m|n)$ whose entries are smooth \textup{(}or holomorphic\textup{)} functions of the even time $t$ with values in $\cC$. Let $X_1(t)$, $X_2(t)$, $X_3(t)$, $X_4(t)$ be four even $(m|n)\times(m|n)$ supermatrices satisfying the \textit{super matrix Riccati equation}
\be\label{Ricc}
\dot X=\beta+\alpha X-X\delta-X\gamma X,
\ee
and such that the differences in Kantor's cross ratio \eqref{(CR)} are invertible at $t=0$. Then, the rational function
\[
\lambda\mapsto\Ber\big((X_1,X_2,X_3,X_4)-\lambda E\big)
\]
does not depend on $t$, so its zeros and poles are first integrals of the system of four equations \eqref{Ricc}; for $n=0$, the coefficients of $\det((X_1,X_2,X_3,X_4)-\lambda E)$ are first integrals.
\end{Proposition}

\begin{proof} Let $g(t)=\begin{pmatrix}a&b\\ c&d\end{pmatrix}$ be the solution of the linear equation $\dot g=Hg$ with $g(0)=E$. If $\begin{pmatrix}P\\ Q\end{pmatrix}$ satisfies the same linear equation and $Q$ is invertible, then $X:=PQ^{-1}$ satisfies equation \eqref{Ricc}; hence, by uniqueness, $X_i(t)=(aX_i(0)+b)(cX_i(0)+d)^{-1}$. For the linear fractional transformation $Y=(aX+b)(cX+d)^{-1}$ and two points $X_i$, $X_j$, a~direct computation gives
\[
Y_i-Y_j=N_i(X_i-X_j)M_j^{-1},\quad\text{where $N_i:=a-Y_ic$ and $M_j:=cX_j+d$;}
\]
it uses only the associativity of the multiplication of supermatrices. Applying this formula to the four differences in the cross ratio \eqref{(CR)} and cancelling, we get
\[
(Y_1,Y_2,Y_3,Y_4)=N_1\,(X_1,X_2,X_3,X_4)\,N_1^{-1}.
\]
The Berezinian is multiplicative, and so the conjugation by $N_1$ does not change $\Ber((Y_1,Y_2,Y_3,Y_4)-\lambda E)$.
\end{proof}

\sssbegin{Remark}[Riccati equations on the other superGrassmannians]\label{RicRem} 1) For $\cC=\Cee$ and $m=1$, $n=0$, Proposition~$\ref{Ric}$ is the classical fact: the cross ratio of four solutions of a~scalar Riccati equation is constant.

2) If $H(t)$ belongs to a~subalgebra of $\fgl(2m|2n)$, e.g., $\fosp$ or $\fpe$, preserving a~subsuperGrassmannian of ${\tt Gr}_{m|n}^{2m|2n}$, then the solutions starting on that subsuperGrassmannian stay on it, and Proposition~$\ref{Ric}$ applies verbatim; this is the super version of the Lagrangian case considered in \cite{Ze}.
\end{Remark}

\ssec{The cross ratio on the superquadrics $Q_{m-2, n}$ and on $E_\alpha$}\label{CRans} Let $J$ be a~Jordan superalgebra with unit $e$, and $\cC$ a~supercommutative superalgebra. The even part of $J\otimes\cC$ is a~Jordan algebra, so the powers of each of its elements $a$ are defined. The element $a$ itself is a~sum of elements of $J_\ev\otimes\cC_\ev$ and of $J_\od\otimes\cC_\od$. If $a$ is in general position, then its powers satisfy a~unique relation of the least degree
\be\label{minpoly}
a^r+c_1a^{r-1}+\dots+c_{r-1}a+c_re=0
\ee
whose \textbf{scalar} coefficients $c_1$, \dots, $c_r$ belong to $\cC_\ev$; the polynomial \[
p_a(\lambda):=\lambda^r+c_1\lambda^{r-1}+\dots+c_r
\]
is the \textit{generic minimal polynomial} of $a$, and its degree $r$, which is the same for all elements in general position, is the \textit{rank} of $J$.

For the superquadrics and for $E_\alpha$, the Jordan superalgebra $J=\fg_{-1}$ is of rank~$2$, and this is what makes the cross ratio computable.

\sssec{The superquadric} Let $\fg=\fosp(m|2n)$ with $m\geq3$ preserve the form with Gram matrix $\antidiag(1,1)\oplus\Gamma$ in a~basis $f_1$, $f_2$, $f_3,\dots$, where $\Gamma$ is the Gram matrix of an even non-degenerate symmetric bilinear form on a~superspace of $\sdim=(m-2|2n)$. The $\Zee$-grading of $\fg$ given by the eigenvalues of $\ad_d$, where $d:=\diag(1,-1,0_{m-2|2n})$, is of depth $1$, and 
\[
\fg_0=\mathfrak{cosp}(m-2|2n):=\fosp(m-2|2n)\oplus\Cee d. 
\]
The Jordan superalgebra $\fg_{-1}$ is the \textit{spin factor} $J=\Cee e\oplus V$, where $V$ is a~superspace of ${\sdim V=(m-3|2n)}$ with an even non-degenerate symmetric bilinear form $(\cdot,\cdot)$, and 
\[
\text{$e\circ x=x$, \ \ $v\circ w=(v,w)e$ for all $v,w\in V$.}
\]
The name \lq\lq spin factor\rq\rq\ comes from the theory of Jordan operator algebras, see \cite[Ch.~6]{HOS}. Let $H$ be a~real Hilbert space and $a$ a~linear map of $H$ into a~$C^*$-algebra satisfying the canonical anticommutation relations of fermions. Then the self-adjoint elements $b(f):=a(f)+a(f)^*$ satisfy $b(f)\circ b(g)=(f,g)1$, so the images of an orthonormal basis of $H$ are symmetries \textup{(}their squares are $1$\textup{)} that pairwise anticommute; such a~set is called a~\textit{spin system}, the Pauli spin matrices being the spin system of cardinality $3$. The Jordan algebra $\Ree1\oplus b(H)$ generated by a~spin system is called a~\textit{spin factor}, a~\lq\lq factor\rq\rq\ because its center is $\Ree1$, as for von Neumann's factors. The infinite-dimensional spin factors were first studied by Topping \cite{To1}, whose next paper \cite{To2} carries the term in its title, and by St{\o}rmer \cite{St}; the finite-dimensional ones constitute the only infinite series of formally real Jordan algebras, other than the matrix ones, in the classification by Jordan, von Neumann and Wigner \cite{JvNW}. McCrimmon \cite{Mcb} also calls them \textit{quadratic factors}.
For every even element $x=x_0e+v$ of $J\otimes\cC$, set
\be\label{normQ}
n(x):=x_0^2-(v,v),\qquad \bar x:=x_0e-v.
\ee
Recall that an even element $x$ of a~unital Jordan superalgebra $J\otimes\cC$ is \textit{invertible} if there is an even $y$ such that $x\circ y=e$ and $(x\circ x)\circ y=x$; this $y$ is unique and is denoted $x^{-1}$, see \cite{Mcb}. Here, $x\circ x-2x_0x+n(x)e=0$, so $x$ is invertible if and only if $n(x)$ is, and then $x^{-1}=n(x)^{-1}\bar x$. Hence, the super version of \textit{Hua's identity}, see \cite{Hua}, and, for Jordan algebras, \cite{Mcb}, holds:
\be\label{Hua}
n(x^{-1}-y^{-1})=\frac{n(x-y)}{n(x)\,n(y)}\quad\text{for all invertible even $x, y$}.
\ee
Indeed, $x^{-1}-y^{-1}=\frac{n(y)\bar x-n(x)\bar y}{n(x)n(y)}$, and
\be\label{Huapf}
n\bigl(n(y)\bar x-n(x)\bar y\bigr)=n(x)n(y)\bigl(n(x)+n(y)-2\langle x, y\rangle\bigr)=n(x)n(y)n(x-y),
\ee
where $\langle x, y\rangle:=x_0y_0-(v, w)$ for $y=y_0e+w$ is the polarization of $n$. The only property of the form used is $(v, w)=(w, v)$ for even $v, w\in V\otimes\cC$, i.e., the symmetry of $(\cdot,\cdot)$ in the super sense.

\sssbegin{Theorem}\label{ThQ} Let $m\geq3$, and let $A, B, C, D\in(J\otimes\cC)_\ev$ be four points of $Q_{m-2, n}$ in general position. Set
\be\label{uv}
u:=\frac{n(A-B)\,n(C-D)}{n(C-B)\,n(A-D)},\qquad v:=\frac{n(A-C)\,n(D-B)}{n(C-B)\,n(A-D)}.
\ee
Then, $u$ and $v$ are invariant under $\OSp(m|2n)$, and every invariant rational function of four points in general position is a~rational function of $u$ and $v$. \end{Theorem}

\begin{proof} On the set of $\cC$-points, the supergroup is generated by the translations $x\mapsto x+a$, by the structure supergroup with Lie superalgebra $\fg_0=\mathfrak{cosp}(m-2|2n)$, acting linearly and multiplying $n$ by a~character, and by the inversion $x\mapsto -x^{-1}$. The translations do not change differences. Each point occurs in the numerator and in the denominator of $u$ and of $v$ equally many times, so the character cancels and, by Hua's identity \eqref{Hua}, so do the factors $n(A)$, $n(B)$, $n(C)$, $n(D)$ produced by the inversion.

The transformation $x\mapsto(x-B)^{-1}-(D-B)^{-1}$ of the supergroup sends $B$ to infinity, $D$ to $0$, and $A$, $C$ to
\be\label{acprime}
a':=(A-B)^{-1}-(D-B)^{-1},\qquad c':=(C-B)^{-1}-(D-B)^{-1}.
\ee
Applying Hua's identity \eqref{Hua} to $x=A-B$, $C-B$, $D-B$ gives
\be\label{uvprime}
u=\frac{n(c')}{n(a')},\qquad v=\frac{n(a'-c')}{n(a')}.
\ee
The stabilizer of the pair (infinity, $0$) is the structure supergroup, whose orbit through $e$ is open; let $g$ in it send $c'$ to $e$. The stabilizer of $e$ in the structure supergroup is the automorphism supergroup of $J$, i.e., the orthosymplectic supergroup of $V$, so the invariants of four points are the invariants of the orthosymplectic supergroup on the remaining point $ga'=a_0e+w$. These are generated by $a_0$ and $(w, w)$, see \cite{Ser}. Since the character of $g$ cancels in formulas \eqref{uvprime}, we have
\be\label{uvnorm}
u^{-1}=n(ga')=a_0^2-(w,w),\qquad vu^{-1}=n(ga'-e)=(a_0-1)^2-(w,w),
\ee
whence $a_0=\frac12(1+u^{-1}-vu^{-1})$ and $(w,w)$ are rational functions of $u$ and $v$.\end{proof}

\sssbegin{Remark}[The invariants $u$, $v$ and Kantor's characteristic polynomial]\label{RemQ} 1) The two invariants \eqref{uv} are the values at $\lambda=0$ and $\lambda=1$ of Kantor's characteristic polynomial \eqref{CRd}. Indeed, for the cross ratio \eqref{(CR)} of four square supermatrices, set $X:=A-B$, $Y:=C-B$, $W:=D-B$; then $A-D=X-W$ and $C-D=Y-W$, and
\be\label{E-M}
(A-D)-XY^{-1}(C-D)=(XY^{-1}-E)W=(A-C)(C-B)^{-1}(D-B),
\ee
so
\be\label{E-M2}
E-(A, B, C, D)=(A-C)(C-B)^{-1}(D-B)(A-D)^{-1}.
\ee
Hence, $\det(A, B, C, D)$ and $\det(E-(A, B, C, D))$ are given by formulas \eqref{uv} with $\det$ in place of $n$, and so are the Berezinians if $mn\neq0$. On a~quadric, the rank of the Jordan superalgebra $J$ is equal to $2$ and nothing else is left: the cross ratio is the quadratic polynomial $\lambda^2-(1+u-v)\lambda+u$, and its roots $z$ and $\bar z$ are the analogs of the eigenvalues of the matrix \eqref{(CR)}.

2) For $n=0$ and $m-2\geq2$, the pair $u=z\bar z$, $v=(1-z)(1-\bar z)$ is the pair of conformal cross ratios on which the four-point functions of conformal field theory depend; for $m-2=2$, the quadric is $\Cee\Pee^1\times\Cee\Pee^1$, and $z$ and $\bar z$ are the classical cross ratios of the projections of the four points onto the two factors. The invariant $u$ is the ratio of the values of the quadratic form $(x,x):=n(x)$ of $J$
\be\label{CRQ}
u=\frac{(A-B, A-B)}{(C-B, C-B)}\frac{(C-D, C-D)}{(A-D, A-D)},
\ee
which suggests itself as the cross ratio on a~quadric; Theorem~$\ref{ThQ}$ says that $v$ is the only other invariant, in every superdimension.

3) For $m=3$, the space $V$ is purely odd, so the invariant $(w, w)$ is nilpotent, and over the underlying $\Cee\Pee^1$ the pair $u, v$ reduces to a~single classical cross ratio.\end{Remark}

\sssec{The exceptional superdomain $E_\alpha$} Let $\fg=\fosp(4|2;\alpha)$ be realized with the Cartan matrix
\be\label{CMosp42}
\begin{pmatrix}2&-1&0\\ -1&0&-\alpha\\ 0&-1&2\end{pmatrix},
\ee
and let the Chevalley generators $X_1^\pm$ be of degree $\pm1$ and $X_2^\pm$, $X_3^\pm$ of degree $0$. The coefficient of $\alpha_1$ in the highest root $\alpha_1+2\alpha_2+\alpha_3$ is equal to $1$, so this $\Zee$-grading is of depth $1$; its component $\fg_{-1}$ is spanned by the root vectors of the roots $-\alpha_1$, $-\alpha_1-\alpha_2$, $-\alpha_1-\alpha_2-\alpha_3$, $-\alpha_1-2\alpha_2-\alpha_3$, and $\fg_0\simeq\fgl(1|2)\simeq\mathfrak{cosp}(2|2)$ is generated by $X_2^\pm$, $X_3^\pm$ and the maximal torus. The Jordan superalgebra structure on $\fg_{-1}$ is $D_t$ of Kac's classification \cite{Ka} with $t$ determined by $\alpha$, see \cite{BC}; $\sdim D_t=(2|2)$, with even basis elements $e_1$, $e_2$, odd ones $x$, $y$, and
\be\label{Dt}
e_i\circ e_j=\delta_{ij}e_i,\quad e_i\circ x=\tfrac12x,\quad e_i\circ y=\tfrac12y,\quad x\circ y=-y\circ x=e_1+te_2,\quad x\circ x=y\circ y=0.
\ee
The unit is $e=e_1+e_2$, and $\fder D_t\simeq\fosp(1|2)$: its even part $\fsl(2)$ acts on $\Span(x, y)$ and kills $e_1$, $e_2$, and its odd part is spanned by
\be\label{DerDt}
\begin{array}{llll}
D_x\colon & e_1\mapsto -y,\ e_2\mapsto y,& x\mapsto -2e_1+2te_2,& y\mapsto0,\\
D_y\colon & e_1\mapsto -x,\ e_2\mapsto x,& y\mapsto 2e_1-2te_2,& x\mapsto0.
\end{array}
\ee

\sssbegin{Theorem}\label{ThE} Let $a=a_1e_1+a_2e_2+\xi x+\eta y$ be an even element of $D_t\otimes\cC$ with $a_1-a_2$ invertible. The rational functions invariant under $\fder D_t\simeq\fosp(1|2)$ are the rational functions of
\be\label{lambdas}
\lambda_1(a):=a_1+\frac{2\xi\eta}{a_2-a_1},\qquad \lambda_2(a):=a_2-\frac{2t\,\xi\eta}{a_2-a_1}.
\ee
Accordingly, the cross ratio of four points of $E_\alpha$ in general position is the pair $\bigl(\lambda_1(ga'),\lambda_2(ga')\bigr)$, where $a'$ and $c'$ are given by formulas \eqref{acprime} with $x^{-1}$ the inverse in $D_t\otimes\cC$ and $g$ is an element of the structure supergroup sending $c'$ to $e$; it is defined up to the automorphisms of $D_t$ not in the identity component. \end{Theorem}

\begin{proof} An invariant of the even part of $\fosp(1|2)$ is of the form 
\[
{F=f(a_1, a_2)+\xi\eta\, g(a_1, a_2)}. 
\]
Under $\eps D_x$ with an odd parameter $\eps$, we have
\be\label{epsDx}
a_1\mapsto a_1+2\eps\xi,\quad a_2\mapsto a_2-2t\eps\xi,\quad \xi\mapsto\xi,\quad \eta\mapsto\eta+\eps(a_2-a_1),
\ee
so the invariance of $F$ reads $\eps\xi\bigl(2\partial_1f-2t\partial_2f-(a_2-a_1)g\bigr)=0$, and $D_y$ gives the same condition. Hence, $g=\frac{2(\partial_1f-t\partial_2f)}{a_2-a_1}$, i.e., $F=f(\lambda_1,\lambda_2)$. The reduction of four points $A$, $B$, $C$, $D$ of $E_\alpha$ in general position to the one point $ga'$ is that of the proof of Theorem~$\ref{ThQ}$, since the structure superalgebra $\fg_0=\mathfrak{cosp}(2|2)$ has an open orbit through $e$ with stabilizer $\fder D_t$.\end{proof}

Over the underlying $\Cee\Pee^1\times\Cee\Pee^1$, the pair $(\lambda_1,\lambda_2)$ reduces to $(a_1, a_2)$, the two classical cross ratios of the projections onto the factors, exactly as for the quadric $Q_{2, n}$; the odd coordinates only add the nilpotent corrections in formulas \eqref{lambdas}, and there is no odd invariant.

\ssec{What replaces Hua's formula}\label{CRgen} In Subsec.~\ref{CRans}, Hua's identity \eqref{Hua} serves only to express $u$ and $v$ through the four points themselves; the invariants are found at one normalized point, and this normalization does not use the rank of $J$ defined in Subsec.~\ref{CRans}. Let $J$ be a~finite-dimensional Jordan superalgebra with unit $e$ whose structure superalgebra~$\fg_0$ has an open orbit through $e$ with stabilizer $\fder J$; this is the case for $Q_{m-2, n}$, $E_\alpha$ and $AB(3)$. The transformation of the proof of Theorem~$\ref{ThQ}$ sends four points in general position to $\infty$, $0$, $e$ and $ga'$, so the invariants of four points are the $\fder J$-invariants of the one even element $ga'\in J\otimes\cC$.

Let $a$ be an even element of $J\otimes\cC$ in general position, and $p_a$ its generic minimal polynomial, of degree $r$, see formula \eqref{minpoly}. The roots of $p_a$ over the underlying algebra are simple, so they lift uniquely to roots in $\cC_\ev$; we call them the \textit{lifted roots} of $a$. They are the super analogs of the eigenvalues of Kantor's cross ratio, the matrix \eqref{(CR)}. The derivations of $J$ kill $e$ and commute with the powers of $a$, so the coefficients of the generic minimal polynomial are invariant under $\fder J$, and so are the lifted roots up to the permutations that the identity component of $\Aut J$ can perform.

\sssbegin{Lemma}\label{LemGen} Let the odd part of $\fder J$ map onto $J_\od$ under $D\mapsto D(a_0)$ at a~generic point $a_0\in J_\ev$. Then, a~$\fder J$-invariant function on $(J\otimes\cC)_\ev$ is determined by its restriction to $J_\ev$. In particular, if the invariants of $(\fder J)_\ev$ on $J_\ev$ are functions of the roots of the minimal polynomial, then every $\fder J$-invariant is a~function of the lifted roots.\end{Lemma}

\begin{proof} Let $F$ be invariant and vanish on $J_\ev$. Expand $F=\sum_{k}F_k$, where $F_k$ is of degree $k$ in the odd coordinates, and let $F_k$ be the lowest non-zero term; then $k\geq1$. Under $\eps D$ with $D$ odd and $\eps$ an odd parameter, the part of the variation of $F$ which is of degree $k-1$ in the odd coordinates comes from $F_k$ alone, and it is the derivative of $F_k$ along the odd vector $\eps D(a_0)$. By the hypothesis, these vectors span all odd directions, so $F_k$ does not depend on the odd coordinates, which contradicts $k\geq1$. \end{proof}

The hypothesis of Lemma~$\ref{LemGen}$ holds for $Q_{m-2, n}$ with $m\geq4$, for $D_t$ and for the Jordan superalgebra of $AB(3)$ considered in Subsec.~\ref{CRab}; it fails for $m=3$, where the nilpotent invariant $(w, w)$ of Remark~$\ref{RemQ}$, 3) appears. Theorems $\ref{ThQ}$ and~$\ref{ThE}$ are the cases $r=2$ of the following statement: \textit{the cross ratio of four points is the collection of the lifted roots of $ga'$}, those that the identity component permutes being taken together.

\sssec{The super norm} In the non-super case, the product of the roots is the generic norm $N$, and Hua's identity 
\[
N(x^{-1}-y^{-1})=(-1)^r\frac{N(x-y)}{N(x)N(y)}
\]
turns two values of the characteristic polynomial into ratios of norms of differences, as in formula \eqref{uv}. In the super case, the product of the lifted roots does not satisfy Hua's identity, already for $D_t$. What satisfies it is a~product of \textbf{powers} of the lifted roots, the \textit{super norm}:
\be\label{superN}
N=z\bar z=n\ \text{ for }Q_{m-2, n},\quad N=\lambda_1^{\,t}\lambda_2\ \text{ for }D_t,\quad N=\mu^2(z\bar z)^{-3}\ \text{ for }AB(3),
\ee
where $z$, $\bar z$ are the lifted roots of Remark~$\ref{RemQ}$, $\lambda_1$, $\lambda_2$ are given by formulas \eqref{lambdas}, and $\mu$, $z$, $\bar z$ are the lifted roots of Theorem~$\ref{ThAB}$ below. For $t\notin\Zee$, the function $N$ is many-valued, but only the ratios \eqref{uvN} of its values occur. In each of the three cases, $N$ is a~relative invariant of the structure supergroup and satisfies Hua's identity up to a~sign. Hence, the pair
\be\label{uvN}
u_N:=\frac{N(A-B)\,N(C-D)}{N(C-B)\,N(A-D)},\qquad v_N:=\frac{N(A-C)\,N(D-B)}{N(C-B)\,N(A-D)}
\ee
is invariant, by the proof of Theorem~$\ref{ThQ}$ verbatim, and so is the \textit{relative characteristic function} of the pencil $a'-\lambda c'$:
\be\label{Phi}
\Phi(\lambda):=\frac{N(a'-\lambda c')}{N(c')}=N(ga'-\lambda e)=\prod_{i}(\lambda_i-\lambda)^{k_i},
\ee
where the $\lambda_i$ are the lifted roots of $ga'$ and the $k_i$ are the exponents in formulas \eqref{superN}. Its two values given by Hua's identity are
\be\label{Phi01}
\Phi(0)=u_N^{-1},\qquad \Phi(1)=v_Nu_N^{-1}.
\ee
\textbf{This is what replaces Hua's formula}: the function $\Phi$ of the pencil, of which Hua's identity yields exactly two values through the four points. In rank $2$, these two values determine the lifted roots; for $E_\alpha$, they are $\lambda_1^t\lambda_2$ and $(1-\lambda_1)^t(1-\lambda_2)$ at the point $ga'$, and over the underlying $\Cee\Pee^1\times\Cee\Pee^1$ they are the $t$-th power of the classical cross ratio of the projections onto the first factor times the classical cross ratio of the projections onto the second factor. In rank $3$, the two values do not suffice: the third invariant of $AB(3)$ is not a~ratio of norms of differences, and only the pencil~\eqref{Phi} gives this third invariant through the four points.

\ssec{The cross ratio on $AB(3)$}\label{CRab} Let $\fg=\fab(3)$, of $\sdim=(24|16)$, be realized with the Cartan matrix
\be\label{CMab3}
\begin{pmatrix}2&-1&0&0\\ -3&0&1&0\\ 0&-1&2&-2\\ 0&0&-1&2\end{pmatrix},
\ee
and let the Chevalley generators $X_1^\pm$ be of degree $\pm1$ and $X_2^\pm$, $X_3^\pm$, $X_4^\pm$ of degree $0$. In the realization of the roots by $\eps_1$, $\eps_2$, $\eps_3$, $\delta$ with $(\eps_i,\eps_j)=\delta_{ij}$, $(\delta,\delta)=-3$, the simple roots are $\alpha_1=-\delta$, $\alpha_2=\frac12(\delta-\eps_1-\eps_2-\eps_3)$, $\alpha_3=\eps_3$, $\alpha_4=\eps_2-\eps_3$, and the degree of a~root is its value under the linear functional equal to $-1$ on $\eps_1$ and $\delta$ and to $0$ on $\eps_2$ and $\eps_3$. The values of this functional on the roots are $-1$, $0$, $1$, so this $\Zee$-grading is of depth $1$; $\fg_0\simeq\mathfrak{cosp}(2|4)=\fosp(2|4)\oplus\Cee z$, of $\sdim=(12|8)$, is generated by $X_2^\pm$, $X_3^\pm$, $X_4^\pm$ and the maximal torus, and $\fg_{-1}$ is spanned by the root vectors of the roots $\eps_1$, $\eps_1\pm\eps_2$, $\eps_1\pm\eps_3$, $\delta$ (even) and $\frac12(\delta+\eps_1\pm\eps_2\pm\eps_3)$ (odd). The Jordan superalgebra $\fg_{-1}$ is Kac's exceptional $K_{10}$ of $\sdim=(6|4)$, see \cite{Ka}. Following \cite{BE}, realize it as $K_{10}=\Cee1\oplus K\otimes K$, where $K=\Cee e\oplus\Cee x\oplus\Cee y$ is Kaplansky's superalgebra, $p(x)=p(y)=\od$:
\be\label{K3}
e\circ e=e,\quad e\circ x=\tfrac12x,\quad e\circ y=\tfrac12y,\quad x\circ y=-y\circ x=\tfrac12e,\quad x\circ x=y\circ y=0,
\ee
with the supersymmetric form $(e|e)=\frac12$, $(x|y)=-(y|x)=\frac12$, the element $1$ being the unit, and
\be\label{K10}
(a\otimes b)\circ(c\otimes d):=(-1)^{p(b)p(c)}\bigl(a\circ c\otimes b\circ d-\tfrac34(a|c)(b|d)\,1\bigr).
\ee
The derivations of $K$ extend to $K_{10}$ by $D(a\otimes b):=D(a)\otimes b$ and $D'(a\otimes b):=(-1)^{p(D')p(a)}a\otimes D'(b)$, and
\be\label{DerK10}
\fder K\simeq\fosp(1|2),\qquad \fder K_{10}=\fder K\otimes1\oplus1\otimes\fder K\simeq\fosp(1|2)\oplus\fosp(1|2);
\ee
the superalgebra $\fder K_{10}\oplus L_{K_{10}}$ is of $\sdim=(12|8)=\sdim\mathfrak{cosp}(2|4)$, where $L_{K_{10}}$ is the space of the operators of multiplication. Write an even element of $K_{10}\otimes\cC$ as
\be\label{aK10}
a=a_0\,1+a_1\,e\otimes e+\sum_{i, j\leq2}m_{ij}\,x_i\otimes x_j+\xi_1\,e\otimes x+\xi_2\,e\otimes y+\xi_3\,x\otimes e+\xi_4\,y\otimes e,
\ee
where $x_1:=x$, $x_2:=y$, and set $M:=(m_{ij})\in\Mat(2;\cC_\ev)$,
\be\label{DeltaOmega}
\begin{array}{l}
\Delta:=a_1^2-\det M,\qquad \Theta:=\xi_1\xi_2+\xi_3\xi_4,\\
\Omega:=a_1\Theta+m_{22}\xi_1\xi_3-m_{21}\xi_2\xi_3-m_{12}\xi_1\xi_4+m_{11}\xi_2\xi_4.
\end{array}
\ee

\sssbegin{Theorem}\label{ThAB} Let $a$ be an even element \eqref{aK10} of $K_{10}\otimes\cC$ with $a_1$ and $\Delta$ invertible. Its generic minimal polynomial is of degree $3$, and its lifted roots are
\be\label{rootsK10}
\mu=a_0+\tfrac32\,\delta,\qquad z, \bar z=a_0+\tfrac12\,\delta\pm\sqrt{\rho},
\ee
where
\be\label{deltarho}
\begin{array}{l}
\delta:=\frac12a_1-\frac{\Omega}{2\Delta}-\frac{3a_1}{2\Delta^2}\,\xi_1\xi_2\xi_3\xi_4,\\[2mm]
\rho:=\frac14\det M+\frac14\Theta-\frac{a_1\Omega}{4\Delta}-\frac{2a_1^2+\det M}{4\Delta^2}\,\xi_1\xi_2\xi_3\xi_4.
\end{array}
\ee
The rational functions invariant under $\fder K_{10}$ are the rational functions of $a_0$, $\delta$ and $\rho$. Accordingly, the cross ratio of four points of $AB(3)$ in general position is the triple $\bigl(a_0(ga'),\delta(ga'),\rho(ga')\bigr)$, where $a'$ and $c'$ are given by formulas \eqref{acprime} with $x^{-1}$ the inverse in $K_{10}\otimes\cC$ and $g$ is an element of the structure supergroup sending $c'$ to $1$. There is no odd invariant. \end{Theorem}

\begin{proof} The minimal polynomial and its lifted roots are obtained by solving 
\[
\text{$a^3-s_1a^2+s_2a-s_3\,1=0$ for $s_1$, $s_2$, $s_3\in\cC_\ev$.}
\]
The even part of $\fder K_{10}\simeq\fosp(1|2)\oplus\fosp(1|2)$ is $\fsl(2)\oplus\fsl(2)$ acting on $M$ by 
\[
\text{$M\mapsto XM+MY^{t}$ for $(X,Y)\in\fsl(2)\oplus\fsl(2)$};
\]
it kills $1$ and $e\otimes e$; its invariants on the even part of $K_{10}$ are the functions of $a_0$, $a_1$ and $\det M$, i.e., of the underlying $\mu$, $z+\bar z$ and $z\bar z$. At a~generic even point, the odd part of $\fder K_{10}$, of dimension $4$, maps onto the odd part of $K_{10}$, of dimension $4$, so Lemma~$\ref{LemGen}$ applies. The reduction of four points to one point is that of the proof of Theorem~$\ref{ThQ}$. \end{proof}

Over the underlying $\Cee\Pee^1\times Q_5$, the triple $(a_0,\delta,\rho)$ of Theorem~$\ref{ThAB}$ reduces to the three roots of the underlying cubic, which factors as 
\[
(\lambda-a_0-\tfrac34a_1)\bigl((\lambda-a_0-\tfrac14a_1)^2-\tfrac14\det M\bigr): 
\]
the linear factor gives the classical cross ratio of the projections onto $\Cee\Pee^1$, and the quadratic factor gives the pair of Theorem~$\ref{ThQ}$ for $Q_5$. The odd coordinates only add the nilpotent corrections in formulas \eqref{deltarho}. 

\ssec{The periplectic quadric}\label{CRpe} Let $\omega$ be the odd symmetric form on $\Cee^{n|n}$ given, for $v=(x|\theta)$ and $v'=(x'|\theta')$, by $\omega(v, v'):=\sum_{i\leq n}(x_i\theta'_i+\theta_ix'_i)$, and let
\be\label{pe}
\fpe(n)=\Bigl\{\begin{pmatrix}A&B\\ C&-A^t\end{pmatrix}\Bigm| B=B^t,\ C=-C^t\Bigr\},\quad \fspe(n)=\{X\in\fpe(n)\mid \tr A=0\}.
\ee
The superdomain $PeQ_{n-1}$ of $\sdim=(n-1|n-1)$ consists of the even lines $\Cee v$ isotropic with respect to $\omega$, i.e., with $\sum x_i\theta_i=0$. Its underlying manifold is $\Cee\Pee^{n-1}$.

\sssbegin{Theorem}\label{ThPe} \textup{1)} For every $k$, every $\mathrm{Pe}(n)$-invariant function of $k$ points of $PeQ_{n-1}$ is a~constant.

\textup{2)} For every $k$, every $\mathrm{SPe}(n)$-invariant function of $k$ points of $PeQ_{n-1}$ is a~constant plus a~nilpotent function; in particular, no invariant of four points of $PeQ_1$ reduces to the classical cross ratio on the underlying $\Cee\Pee^1$.

\textup{3)} For $n=2$ and four points $v_i=(x_i|\theta_i)$ in general position, set
\be\label{wij}
w_{ij}:=\frac{\omega(v_i, v_j)}{\det(x_i, x_j)}.
\ee
The $w_{ij}$ are odd and $\mathrm{SPe}(2)$-invariant, $w_{ij}+w_{jk}=w_{ik}$, and the $\mathrm{SPe}(2)$-invariant functions of the four points are the linear combinations
\be\label{invPe}
c_0+\sum_{j}c_j\,w_{1j}+\sum_{j<l}c_{jl}\,w_{1j}w_{1l}+f(\kappa)\,w_{12}w_{13}w_{14}
\ee
with constant $c_0$, $c_j$, $c_{jl}$, where $j, l\in\{2, 3, 4\}$, $\kappa$ is the classical cross ratio of the underlying points and $f$ is an arbitrary function. Thus, the cross ratio on $PeQ_1$ is the triple of odd invariants $w_{12}$, $w_{13}$, $w_{14}$, and the classical cross ratio survives only as the coefficient of their product. \end{Theorem}

\begin{proof} At a~point $(x|\theta)$, an element of $\fpe(n)$ acts by $\delta x=Ax+B\theta$ and $\delta\theta=Cx-A^t\theta$. Normalize the representative of a~line by its even part; then the coordinates of a~point are $[x]\in\Cee\Pee^{n-1}$ and the odd coordinates of $\theta$ in the directions orthogonal to $x$.

1) The element of $\mathrm{Pe}(n)$ with $A=a\,1_n$ and $B=C=0$ fixes every point of $\Cee\Pee^{n-1}$ and multiplies the odd coordinates of every point by $a^{-2}$, so an invariant does not depend on the odd coordinates. For a~function $F$ of the underlying points alone, the odd elements with $A=C=0$ give $\delta x_i=\eps B\theta_i$, so the invariance requires the derivative of $F$ in $x_i$ along $B\theta_i$ to vanish for every symmetric $B$; since $\{B\theta\mid B=B^t\}=\Cee^n$ for $\theta\neq0$, the function $F$ is constant.

2) The same argument, applied to the part of an invariant of degree $0$ in the odd coordinates, shows that this part is a~constant; here, the element with $A=a\,1_n$ is not available.

3) Take $x=(1, y)$ and $\theta=s(-y, 1)$, where $s$ is odd; then $w_{ij}=s_i-s_j$. The elements of $\fsl(2)$ act on the $y_i$ by the vector fields $\varphi(y_i)\partial_{y_i}$ with $\varphi$ quadratic and fix every $s_i$; the elements with $B\neq0$ fix every $s_i$ and act by $\delta y_i=\eps s_i\varphi(y_i)$, where $\varphi$ runs over all quadratic polynomials; the elements with $C\neq0$ act by $\delta s_i=\eps c$, the same for all $i$. The invariance under the elements with $C\neq0$ says that $F$ depends on the $s_i$ through the differences $w_{j}:=s_j-s_1$ only, so 
\[
\text{$F=\sum_S g_S(y)w^S$ with $S\subset\{2, 3, 4\}$.}
\]
In the variation under $B$, write $s_i=s_1+w_i$ (with $w_1:=0$): the terms containing $s_1$ say that every $g_S$ is invariant under the diagonal action of $\fsl(2)$; the other terms say that $\sum_{i\in T}\pm\varphi(y_i)\partial_{y_i}g_{T\setminus i}=0$ for every $T$. Choosing $\varphi$ with zeros at the two points of $\{y_2, y_3, y_4\}$ other than $y_i$, we get $\partial_{y_i}g_S=0$ for $i\not\in S$, $i\geq2$. So $g_S$ is an $\fsl(2)$-invariant function of the $|S|+1$ points $y_1$ and $y_i$ with $i\in S$: a~constant if $|S|\leq2$, and a~function of $\kappa$ if $|S|=3$. \end{proof}

For $n\geq3$, part 2) leaves room for nilpotent invariants only. At a~generic point of the underlying configuration of four points, the $\mathrm{Pe}(n)$-orbit has codimension $(1|3)$ for $n=2$, $(0|5)$ for $n=3$ and $(0|6)$ for $n\geq4$, whereas by 1) there are no invariants at all: the odd codimension is absorbed by the torus $a\,1_n$ and by the elements with $B\neq0$, which vanish at the underlying points.

\ssec{Problems}  
\textup{1)} What are the analogs of cross ratio 
 related to the curved quadrics, \textup{(}twisted\textup{)} loops, and in stringy cases? Which nilpotent invariants of four points of $PeQ_{n-1}$ does $\mathrm{SPe}(n)$ have for $n\geq3$?

\textup{2)} What is the cross ratio on the classical superspace
associated with the infinite-dimensional exception related with the subalgebra $\fkas^L\subset \fk^L(1|6)$, discovered by Shchepochkina  \cite{Sh14}? (On importance of fast publications: The Lie superalgebras $\fkas^L$, and its version $\fkas$ with polynomial  coefficients, though discovered by Shchepochkina and reported at conferences in 1996, are often called in the literature ``CK" superalgebras.)

\Subsection{Acknowledgements} I am thankful to
J.~Lukierski for hospitality, to 
B.~Zupnik, Ho.~Khudaverdian and
A.~Vaintrob for inspiring discussions, P.~Grozman, V.~Serganova, V.~Shander and I.~Shchepochkina for
help. Financial support of the grant AD 065 NYUAD is acknowledged.

\section{Tables (\cite{S2, LSV})} \label{Tab}

\begin{table*}[!ht]
\caption{Gradings of twisted loop
(super)algebras corresponding to hermitian superdomains}
\label{table:1}
\scriptsize
\setlength{\tabcolsep}{3pt}
\[
\renewcommand{\arraystretch}{1.4}
\begin{tabular}{|c|c|c|c|}
\hline
$\mathfrak{g}_{\varphi}^{(m)}$ &$\varphi$ & grading elements from
$\mathfrak{h}$& $(\mathfrak{g}_{\varphi}^{(m)})_0$ \cr
\hline
$\mathfrak{sl}(2m|2n)^{(2)}$&\twol{$(-st)\circ$}{${\rm Ad}_{{\rm diag}(\Pi
_{2m}, J_{2n})}$} &${\rm diag}(1_m, -1_m, 1_n,
-1_n)$&$\mathfrak{sl}(m|n)^{(1)}$\cr

$\mathfrak{sl}(2m)^{(2)}$&$(-t)\circ {\rm Ad}_{
\Pi_{2m}}$&${\rm diag}(1_m, -1_m)$&$\mathfrak{sl}(m)^{(1)}$\cr 
$\mathfrak{sl}(2n)^{(2)}
$&$(t)\circ {\rm Ad}_{J _{2n}}$&${\rm diag}(1_n,
-1_n)$&$\mathfrak{sl}(n)^{(1)}$\cr

$\mathfrak{sl}(n|n)^{(2)}$&$\Pi$&${\rm diag}(1_p, 0_{n-p}, 1_p,
0_{n-p})$& \twol{$\mathfrak{s}(\mathfrak{gl}(p|p)^{(2)}_{\Pi}\oplus{}$}{$
\mathfrak{gl}(n-p|n-p)^{(2)}_{\Pi})$}\cr

$\mathfrak{sl}(n|n)^{(2)}$&$\Pi\circ (-st)$&${\rm diag}(1_p, -1_{n-p},
-1_p, 1_{n-p})$& \twol{$\mathfrak{s}(\mathfrak{gl}(p|p)^{(2)}_{\Pi\circ
(-st)}\oplus{}$}{$\mathfrak{gl}(n-p|n-p)^{(2)}_{\Pi\circ (-st)})$}\cr
\hline
$\mathfrak{osp}(2m|2n)^{(2)}$&\twol{$\varphi_{m, n}\circ{}$}{${\rm Ad}_{{\rm
diag}(1_{2m-1}, -1, 1_{2n})}$} &${\rm diag} (J_2, 0_{2(m+n-1)}) $&$
(\mathfrak{cosp}(2m-2|2n))^{(1)}_{\varphi_{m-1, n}}$\cr

$\mathfrak{o}(2m)^{(2)}$& &&$(\mathfrak{co}(2m-2))^{(1)}$\cr

\hline
$\mathfrak{psq}(2n)^{(4)}$&$ (-st)\circ \delta_{\sqrt{-1}}$&${\rm diag}(J_{2n},
J_{2n})$& $\mathfrak{psq}(n)^{(2)}_{\delta_{-1}}$\cr
$(\mathfrak{h}'(0|2n))^{(2)}$&$ A$&& $ (\mathfrak{h}'
(0|2n-2)\oplus \Lambda(2n-2))^{(2)}_A $\cr

$\mathfrak{psq}(n)^{(2)}$&$
\delta_{-1}$&${\rm diag}(1_p, 0_{n-p}, 1_p, 0_{n-p})$& \twol{$\mathfrak{ps}
(\mathfrak{q}(p)^{(2)}_{\delta _{-1}}\oplus{}$}{$
\mathfrak{q}(n-p)^{(2)}_{\delta _{-1}})$}\cr
\hline
\end{tabular}
\]

\end{table*}

Let $\mathfrak{s}(\mathfrak{g})$ be the traceless part of
$\mathfrak{g}$ and $\mathfrak{p}(\mathfrak{g})=\mathfrak{g}/{\rm
center}$; let $\mathfrak{g}_{\varphi}^{(m)}$ be the stationary
subalgebra of the loop algebra with values in $\mathfrak{g}$ singled
out by the order $m$ automorphism $\varphi$ of $\mathfrak{g}$. For
$\mathfrak{g}_{\varphi}^{(m)}$ with the $\Zee$-grading of type
$(d=1)$, the last column of Table 1 contains $(\mathfrak{g}_{\varphi}^{(m)})_{0}$; the maps $-\st:X\mapsto -X^{st}$, the map $\Pi:\begin{pmatrix}a& b\cr c&d\end{pmatrix}\mapsto \begin{pmatrix}d& c\cr
b&a\end{pmatrix}$, and $\delta_{x}:\begin{pmatrix}a& b\cr
c&d\end{pmatrix}\mapsto \begin{pmatrix}a& xb\cr xc&d\end{pmatrix}$; the
automorphism $A$ of $\mathfrak{po}$ is defined on monomials
$f(\theta)$ as $\id$ if $\frac{\partial f}{\partial \theta_1}=0$ and $-\id$
otherwise; $\text{irr}(\ldots)$ is any of the two irreducible components;
$LGr$ and $OGr$ stand for the Lagrangian and orthogonal
Grassmannian, respectively; the dual domain is endowed with an
asterisk as a~left superscript.

In Table~2, let $G$ be a~
simple Lie supergroup, $P$ its parabolic subsupergroup corresponding to
several omitted generators of a~Borel subalgebra (description of these
generators can be found in \cite{GLP}). We say that the homogeneous superspace $G/P$ is of \textit{depth} $d$ and \textit{height} $l$ if such are the depth and height of $\mathfrak{g}={\rm Lie}(G)$ in the
$\Zee$-grading compatible with that of ${\rm Lie}(P)$. 

Note that all
superspaces of Table 2 have an hermitian structure (hence, are of
height 1) except for $PeGr$ (no hermitian structure), $PeQ$ (no hermitian structure,
height 2), $CGr_{0, k}^{0, n}$ and $SCGr_{0, k}^{0, n}$ (no hermitian structure,
height $n-k$; for $SCGr_{0, 1}^{0, n}$, the height of the simple $\fsvect'(0|n)$ is $n-2$).

{\emergencystretch=1.5em The curved superquadric $CQ_{m-2, 0}$ has two series of infinite-dimensional ``stringy''
counterparts with $\mathfrak{h}(0| m)$ replaced with the centerless
$N$-extended Neveu-Schwarz algebra $\mathfrak{k}^L(1|N)$ or Ramond
algebra $\mathfrak{k}^M(1|N)$, and one exceptional analog associated with the exceptional stringy Lie superalgebra $\fkas$ discovered by Shchepochkina, see \cite{Sh14, GLS2}.\par}

In Table 2, an interpretation of the
superGrassmannians is given line-wise: the superGrassmannian of
$p|q$-dimensional subsuperspaces in $\Cee^{m|n}$ and same for $n=m$,
$p=q$; superquadric of $1|0$-dimensional isotropic (relative to a
non-degenerate even form) lines in $\Cee^{m|2n}$; ortholagrangian
superGrassmannian; queerGrassmannian; ``odd'' superquadric (relative
to a non-degenerate odd form) of $1|0$-lines in $\Cee^{n|n}$;
odd-lagrangian superGrassmannian; curved superGrassmannian of
$0|k$-dimensional subsupermanifolds in $\Cee^{0|n}$; curved
superquadric; the two finite-dimensional exceptions. 

\begin{table*}[!ht]
\caption{Classical superspaces of depth 1 }
\label{table:2}
\scriptsize
\setlength{\tabcolsep}{3pt}
\[
\renewcommand{\arraystretch}{1.4}
\begin{tabular}{|c|c|c|c|c|}
\hline
$\mathfrak{g} $&$ \mathfrak{g}_0 $&$
\mathfrak{g}_{-1}$&Underlying&Name of the\cr
&&&domain&superdomain\cr
\hline
$\mathfrak{sl}(m|n)$&$ \mathfrak{s}(\mathfrak{gl}(p|q)\oplus
\mathfrak{gl}(m-p|n-q))$&$\id\boxtimes \id^*$& $Gr_p^m\times
Gr_q^n$&$Gr_{p, q}^{m, n}$\cr
$\mathfrak{psl}(m|m)$&$\mathfrak{ps}(\mathfrak{gl}(p| p)\oplus
\mathfrak{gl}(m-p|m-p))$&$ \id\boxtimes \id^*$&$Gr_p^m\times Gr_p^m
$&$Gr_{p, p}^{m, m}$\cr
\hline
$\mathfrak{osp}(m| 2n)$&$ \mathfrak{cosp}(m-2| 2n)$& $\id$& $Q_{m-2}$&
$Q_{m-2, n}$\cr
\hline
$\mathfrak{osp}(2m| 2n)$& $\mathfrak{gl}(m|n)$ & $\Lambda^2(\id)$&
 ${}^*OGr_m\times LGr_n
$&$OLGr_{m, n}$ \cr \hline $\mathfrak{sq}(n)$&$
\mathfrak{s}(\mathfrak{q}(p) \oplus
\mathfrak{q}(n-p))$&$\text{irr}(\id\boxtimes \id^*)$& $Gr_p^n$&$ QGr_p^n$\cr

$\mathfrak{psq}(n)$&$\mathfrak{ps}(\mathfrak{q}(p) \oplus
\mathfrak{q}(n-p))$& &&\cr 
\hline
$\mathfrak{pe}(n)$&$\mathfrak{cpe}(n-1)$&$ \id$& $\Cee P^{n-1}$&
$PeQ_{n-1}$ \cr $\mathfrak{spe}(n)$&$\mathfrak{cspe}(n-1)$&$ \id$ &&\cr

\hline 
$\mathfrak{pe}(n)$&$\mathfrak{gl}(p|n-p)$&$\Pi(S^2 (\id))$&
$Gr_p^n$& $PeGr_p^n$\cr $(\mathfrak{spe}(n))$&$
(\mathfrak{sl}(p|n-p))$& or $\Pi (\Lambda^2 (\id))$&&\cr 
\hline
\hline 
$\mathfrak{vect}(0| n)$&$\mathfrak{vect}(0| k)\oplus
\mathfrak{gl}(n-k; \Lambda (k))$ &$\Lambda (k)\boxtimes \Pi
(\id)$&curved&$CGr_{0, k}^{0, n}$\cr

$\mathfrak{svect}(0|n)$&$\mathfrak{vect}(0| k)\oplus
\mathfrak{sl}(n-k; \Lambda (k))$& $\Pi (\Vol)$ if $k=n-1$ &curved&$SCGr_{0,
k}^{0, n}$\cr 
\hline 
$\mathfrak{h}(0| m)$&$\mathfrak{h}(0|m-2)\oplus
\Lambda (m-2)\cdot z$&$ \Pi (\id)$&curved&$CQ_{m-2, 0}$ \cr

$\mathfrak{h}'(0|m)$&$\mathfrak{h}'(0|m-2)\oplus \Lambda (m-2)\cdot z$
&$ \Pi (\id)$&curved&\cr 
\hline \hline
$\mathfrak{k}^L(1|n+2)$&$\mathfrak{k}^L(1|n)\oplus \cF(1|n)\alpha\beta$&\twol{$\Pi(\Cee[t^{-1}, t, \theta_1, \dots , \theta_n;$}{$ \alpha, \beta])$}&curved&${\tt K{}^LQ}_{1|n}$\cr
 \twol{$\mathfrak{k}^M(1|n+2)$}{for $n>0$}&$\mathfrak{k}^M(1|n)\oplus \cF(1|n)\alpha\beta$&\twol{$\Pi(\Cee[t^{-1}, t, \theta_1, \dots , \sqrt{t} \theta_n;$}{$ \alpha, \beta])$}&curved&${\tt K{}^MQ}_{1|n}$\cr
\hline \hline
$\mathfrak{osp}(4|2;\alpha)$&$\mathfrak{cosp}(2|2)\simeq
\mathfrak{gl}(2|1)$ & $\id$& $\Cee \Pee^1\times \Cee \Pee^1$
&$E_{\alpha}$\cr 
\hline 
$\mathfrak{ab}(3)$&$\mathfrak{cosp}(2|4)$&$
L_{3\varepsilon_{1}}$ & $\Cee \Pee^1\times Q_5$ &$AB(3)$\cr \hline
\end{tabular}
\]\\ 

\end{table*}


\end{document}